%% file: main.tex
\PassOptionsToPackage{dvipsnames,table}{xcolor}
\PassOptionsToPackage{hidelinks}{hyperref}
\PassOptionsToPackage{ruled,linesnumbered}{algorithm2e}
\documentclass{l4dc2027}
\SetAlCapHSkip{0pt}
\RestyleAlgo{ruled}
\usepackage{times}
\usepackage{xcolor}
\usepackage{textcomp}
\usepackage{wrapfig}
\usepackage{mathtools}
\usepackage{float}
\usepackage{physics}
\usepackage{booktabs}
\usepackage{caption}
\newsavebox{\reachertabbox}
\newlength{\tabw}
\newlength{\colgap}
\makeatletter
\newcommand{\newacronym}[3]{%
  \expandafter\def\csname lfourdc@acr@short@#1\endcsname{#2}%
  \expandafter\def\csname lfourdc@acr@long@#1\endcsname{#3}%
}
\newcommand{\gls}[1]{%
  \@ifundefined{lfourdc@acr@used@#1}{%
    \csname lfourdc@acr@long@#1\endcsname\space
    (\csname lfourdc@acr@short@#1\endcsname)%
    \global\@namedef{lfourdc@acr@used@#1}{}%
  }{\csname lfourdc@acr@short@#1\endcsname}%
}
\makeatother

\newacronym{LTI}{LTI}{Linear Time-Invariant}
\newacronym{DeePC}{DeePC}{Data-enabled Predictive Control}

\definecolor{tablebg}{HTML}{F1F8DF}

\input{user_commands.tex}
\begin{document}

\title[Atlas-DeePC]{Atlas-DeePC: Learning Behavioral Atlases for Predictive Control}
\coltauthor{%
  \Name{Dian Jin}, \Name{Zhiyuan Jin},
  \Name{Jeremy Coulson}%
}

\maketitle

\begin{abstract}
    Data-enabled Predictive Control (DeePC) uses linear combinations of recorded input-output trajectories to compute control actions without a parametric model. 
    For nonlinear systems, these linear combinations may not satisfy the system dynamics, and the computational cost grows with dataset size.
    We propose Atlas-DeePC, which approximates the behavior of a dynamical system using a collection of low-dimensional subspaces, termed a behavioral atlas, to enable computationally efficient prediction and control.
    The atlas is constructed offline by clustering the data and fitting a subspace to each cluster. 
    Online, the controller selects the subspace that best matches recent measurements and solves a predictive control problem parameterized by its orthonormal basis. A graph-based reference planner complements the controller by providing intermediate references across operating regions.
    We illustrate Atlas-DeePC on a planar robotic manipulator, where it achieves higher success rates and substantially lower online computation times than the compared controllers based on online data selection. 
    An additional benchmark of an autonomous racing task is also included. The project website provides code and demonstrations
    \footnote{\url{https://dianjin-frederick.github.io/Atlas-DeePC/}}.
\end{abstract}
    

\input{v3/v3}

\end{document}

%% file: user_commands.tex
\newcommand{\mean}[1]{\underset{#1}{\operatorname{mean}}}
\newcommand{\im}{\mathrm{im~}}
\newcommand{\ini}{\mathrm{ini}}
\newcommand{\rmf}{\mathrm{f}}
\newcommand{\Tf}{T_\mathrm{f}}
\newcommand{\Tini}{T_\mathrm{ini}}
\newcommand{\Tsim}{T_\mathrm{sim}}
\newcommand{\lam}{\lambda}
\newcommand{\omg}{\omega}

\newcommand{\R}{\mathbb{R}}
\newcommand{\Z}{\mathbb{Z}}
\newcommand{\Zp}{\mathbb{Z}_{\geq 0}}

\newcommand{\calB}{\mathcal{B}}
\newcommand{\calC}{\mathcal{C}}
\newcommand{\calE}{\mathcal{E}}
\newcommand{\calG}{\mathcal{G}}
\newcommand{\calH}{\mathcal{H}}
\newcommand{\calI}{\mathcal{I}}
\newcommand{\calJ}{\mathcal{J}}
\newcommand{\calP}{\mathcal{P}}
\newcommand{\calR}{\mathcal{R}}
\newcommand{\calS}{\mathcal{S}}
\newcommand{\calT}{\mathcal{T}}
\newcommand{\calU}{\mathcal{U}}
\newcommand{\calV}{\mathcal{V}}
\newcommand{\calY}{\mathcal{Y}}

\DeclareMathOperator*{\diag}{diag}

\ifdefined\theoremstyle
  \theoremstyle{remark}
\fi
\newtheorem{problem}{Problem}

%% file: v3/v3.tex
\section{Introduction}
Controlling nonlinear systems requires handling dynamics that change across operating regions while respecting constraints. Model predictive control addresses constraints explicitly \citep{camacho2007model}, but it relies on an accurate model, which can be difficult to derive. 
Data-driven control methods instead design controllers directly from input-output data, 
bypassing explicit model identification. Among these approaches, 
\gls{DeePC} represents the trajectories of a \gls{LTI} system as a linear subspace spanned by sufficiently informative input-output data \citep{willems2005note}, and uses linear combinations of these trajectories to predict future trajectories in a receding-horizon controller \citep{coulson2019data, berberich2021data}.
However, the input-output trajectories of a nonlinear system generally do not form a linear subspace, and the computational cost grows with the size of the dataset. 

Existing extensions of data-driven predictive control to nonlinear systems follow two broad strategies \citep{berberich2025overview}: exploiting linearity in higher-dimensional coordinates, as in kernel and Koopman formulations \citep{korda2018linear, huang2023robust}, 
or exploiting local linearity of the dynamics, as in predictors built from selected data near the current operating point \citep{berberich2022lineardata, naf2025choose, beerwerth2026less}. 
Global representations require choosing kernels or observables rich enough to capture the nonlinearity across all operating regions, which is difficult without prior knowledge of the system \citep{brunton2016koopman}. 
Local approaches select the data most relevant to the current operating point \citep{naf2025choose, beerwerth2026less}, which requires a search over the full dataset online. 
Alternatively, the data can be partitioned offline into local datasets, either according to a prescribed scheduling variable \citep{guerrero2025gain, zieglmeier2025gain} or
by clustering input-output trajectories into local linear sub-behaviors, as in big data-driven predictive control (BDPC) \citep{han2024big}. 
The former requires system knowledge to select a scheduling variable that captures the dominant nonlinearity \citep[Assumption 2]{zieglmeier2025gain}. 
BDPC approaches construct local sub-behaviors for one-step predictive control \citep{han2024big}, 
with extensions that impose contraction conditions \citep{han2025big} or interpolate sub-behaviors and select 
paths based on uncertainty \citep{han2026big}.
Their one-step formulation \citep[Eq. (34)]{han2026big} limits the controller's prediction horizon and may impair tracking performance. 

We develop Atlas-DeePC, which represents the data by a collection of local low-dimensional subspaces, termed \emph{behavioral charts}, which together form a \emph{behavioral atlas} by analogy with manifold atlases \citep{lee2012manifold}. 
Inspired by manifold learning \citep{roweis2000nonlinear, zhang2004principal} and the local sub-behavior viewpoint \citep{han2024big}, 
we construct the atlas offline by clustering the data and fitting a low-dimensional subspace to each cluster. 
Online, the controller selects the subspace that best matches recent measurements and uses it to formulate and solve a predictive control problem. 
While \citet{eysenbach2019search} plan over stored observations
using learned distances, we plan over behavioral charts
using observed transitions and subspace distances to
generate intermediate tracking references.
Although Atlas-DeePC and BDPC use local behavioral representations, they differ in their online control and reference planning formulations. Appendix \ref{app:related-work} provides a detailed comparison. 
The main contributions of this work are
\begin{enumerate}
    \item[(i)] We propose Atlas-DeePC, a predictive control framework that combines local low-dimensional representations of input-output trajectories of a nonlinear system with online chart selection. 
    \item[(ii)] We develop a graph-guided reference planner that uses transitions observed in offline data to generate intermediate tracking references, guiding the controller through nearby charts toward targets beyond the operating region represented by the current chart.
    \item[(iii)] We validate Atlas-DeePC on planar robotic reaching and autonomous racing tasks. 
    The reaching benchmark demonstrates higher success rates and substantially lower online computation times than the compared data-selection-based predictive controllers.
\end{enumerate}

The rest of the paper is organized as follows. Section 2 formulates the problem and reviews DeePC. 
Section 3 introduces Atlas-DeePC. 
Section 4 presents the planar robotic manipulator benchmark. 
Section 5 introduces graph-guided reference planning. 
Section 6 evaluates reference planning and compares Atlas-DeePC with existing methods.
Section 7 concludes the paper. 
Appendix A extends the discussion of related work.
Appendix B reports additional manipulator experiments, and Appendix C presents the autonomous vehicle example.
\paragraph{Notation.}
Let $\Zp, \Z_{>0}$ denote the sets of nonnegative and positive integers, respectively.
For a signal $s\colon\Zp\to\R^q$ and indices $i,j\in\Zp$ with $i\leq j$, we define $s_{[i,j]}\coloneqq(s_i,s_{i+1},\dots,s_j)\in\R^{q(j-i+1)}$. For sets $X$ and $Y$, let $Y^X$ denote the set of all maps from $X$ to $Y$. 
For a matrix $A=[a_{ij}]\in\R^{r\times s}$ and index sets $\calI\subseteq\{1,\dots,r\}, \calJ \subseteq \{1, \dots, s\}$, $A(\calI,\calJ)$ denotes the submatrix consisting of the 
rows indexed by $\calI$ and the columns indexed by $\calJ$, arranged in increasing index order, and $A(\calI,:)$ and $A(:,\calJ)$ denote row and column submatrices, respectively. If $A$ is a vector, we omit the colon. 
For consecutive indices, $a:b$ denotes the index set $\{a,a+1,\ldots,b\}$. 
Given a finite sequence $s_{[0,M-1]}$ and a depth $1\leq L\leq M$, its block Hankel matrix is
$\calH_L(s_{[0,M-1]})\coloneqq\allowbreak
\bigl[
s_{[0,L-1]}\;\allowbreak
s_{[1,L]}\;\allowbreak
\cdots\;\allowbreak
s_{[M-L,M-1]}
\bigr]
\allowbreak\in\R^{qL\times(M-L+1)}$.
We write $I_n$ for the $n\times n$ identity matrix and $\mathbf{1}_n$ for the vector of $n$ ones. The notation $\operatorname{diag}(A_1,\ldots,A_n)$ denotes
the block-diagonal matrix with scalar or matrix blocks
$A_1,\ldots,A_n$.
\section{Problem Formulation}
\label{sec:problem-formulation}
A \textit{dynamical system} is a 3-tuple $(\Z_{\geq 0}, \R^q, \calB)$, where $\calB \subseteq (\R^q)^{\Z_{\geq 0}}$ is the \textit{behavior} \citep{polderman1998introduction}. Each $w \in \calB$ is called a \textit{trajectory}.
The length-$L$ \textit{restricted behavior} is defined by $\calB_L \coloneqq \{w_{[0, L-1]}: w \in \calB\}$. 
Consider the state-space system
\begin{equation}\label{eq:ss}
    x_{t+1}=f(x_t, u_t), \quad y_t=h(x_t, u_t),
\end{equation}
where $f:\R^n \times \R^m \to \R^n, h: \R^n \times \R^m \to \R^p$ and $x_t \in \R^n, u_t \in \R^m, y_t \in \R^p$ denote the state, input, output, respectively. 
We set $q=m+p$ and define $w_t = (u_t, y_t) \in \R^q$.
The length-$L$ restricted behavior of \eqref{eq:ss} is defined by
\[
\calB_L = \left\{w_{[0, L-1]} \in \R^{(m+p)L} \mid \exists x_{[0, L]} \text{ such that \eqref{eq:ss} holds for all }t \in [0, L-1] \right\}.
\]
In the case that \eqref{eq:ss} is \gls{LTI} of order $n$,  its restricted behavior $\calB_L$ is a linear subspace of $\R^{(m+p)L}$ with $\dim \calB_L = mL+n$ for large enough $L$ \citep{willems2005note}. A basis of $\calB_L$ can be computed from data, without identifying the system parameters \citep{markovsky2022identifiability}.
Let $w_{[0,T-1]}$ be a trajectory of \eqref{eq:ss} and let $H=[w_{[0,L-1]} \cdots w_{[T-L,T-1]}]$,
then $\rank H = mL+n$ if and only if
\(\im H = \calB_L\).
This data-based representation of a behavior underpins \gls{DeePC} \citep{coulson2019data}.
\paragraph{Data-enabled predictive control.} 
Let $\Tini$ and $\Tf$ be positive integers and let
$L=\Tini+\Tf$. Define $P$ to be the permutation matrix such that 
\(
P w_{[0,L-1]} =
\bigl[
u_{[0,\Tini-1]}^\top \allowbreak \,
y_{[0,\Tini-1]}^\top \allowbreak \,
u_{[\Tini,L-1]}^\top \allowbreak \,
y_{[\Tini,L-1]}^\top
\bigr]^\top
\).
At time $t$, the sequences $u_{[t-\Tini,t-1]}$ and
$y_{[t-\Tini,t-1]}$ are the most recent $\Tini$
applied inputs and measured outputs, respectively.
Let $u^\rmf \in \R^{m\Tf}, y^\rmf \in \R^{p\Tf}$ denote the future input and output sequences. 
Let $y^{\mathrm{ref}}\in\R^{p\Tf}$ be an output reference, 
$Q$ be a positive semidefinite matrix and $R$ be a positive definite matrix. The \gls{DeePC} optimal control problem is 
\begin{equation}\label{eq:deepc}
\begin{aligned}
    &\min_{\substack{
        g, \sigma, u^\rmf, y^\rmf}}
    \quad
    \left(y^\rmf-y^{\mathrm{ref}}\right)^\top Q
    \left(y^\rmf-y^{\mathrm{ref}}\right)
    + \left(u^\rmf\right)^\top R u^\rmf + \lam_\ini \|\sigma\|_2^2 \\
    &\mathrm{s.t.}\quad
    PH g
    =
    \begin{bmatrix}
        u_{[t-\Tini,t-1]} \\
        y_{[t-\Tini,t-1]} + \sigma \\
        u^\rmf \\ 
        y^\rmf
    \end{bmatrix}, \quad u^\rmf \in \calU, \quad y^\rmf \in \calY,
\end{aligned}
\end{equation} 
where $\calU, \calY$ are convex constraint sets, $\lam_\ini \geq 0$, and $u^\rmf, y^\rmf$ are decision variables representing the future input-output trajectory prediction.
Let $u_t^{\rmf, *} \in \R^{m\Tf}$ denote the optimal input sequence obtained from solving \eqref{eq:deepc} at time $t$. Applying the first optimal input $u_t \coloneqq u_t^{\rmf, *}(1:m)$ and repeating the optimization at each time step defines the \gls{DeePC} policy $u_t = \pi_H(w_{[t-\Tini,t-1]}, y^\mathrm{ref})$. 
Since the policy $\pi_H$ depends entirely on the subspace $\im H$ \citep{zhang2023dimension, jin2026sensitivity}, we also write $\pi_H$ as $\pi_{\im H}$.
When $\im H = \calB_L$, the trajectory constraint in \eqref{eq:deepc} ensures that $P^\top
\bigl(
u_{[t-\Tini,t-1]},\allowbreak\,
y_{[t-\Tini,t-1]}+\sigma,\allowbreak\,
u^\rmf,\allowbreak\,
y^\rmf
\bigr)
\in\calB_L$.
However, for nonlinear systems, $\calB_L$ need not be a linear subspace.
Consequently, $\im H$ may contain trajectories that cannot be generated by \eqref{eq:ss} and some trajectories of \eqref{eq:ss} may not be generated by $\im H$. 
The manifold viewpoint of nonlinear trajectories \citep{hauser1998trajectory, yan2026behavioral, berberich2022linearmpc, bayer2026trajectory} and trajectory clustering into local sub-behaviors \citep{han2024big, han2025big, han2026big} motivate approximating $\calB_L$ locally by low-dimensional linear subspaces using tools from manifold learning \citep{roweis2000nonlinear, zhang2004principal}. We formulate the following learning and control problem.
\begin{problem}
    Given a collection of trajectories $\{w^{(j)}_{[0,N_j-1]}\}_{j=1}^N$ from \eqref{eq:ss}, we wish to:
    (i) learn $K$ linear subspaces
    $\{\calS_k\}_{k=1}^K$ that locally approximate $\calB_L$;
    (ii) design a selection rule that chooses a subspace
    $\calS_{k_t}$ at each time $t$ based on online measured trajectory data; and
    (iii) design a predictive controller $\pi_{\calS_{k_t}}$
    that tracks an output reference while satisfying input
    and output constraints.
\end{problem}
\section{Learning and control with behavioral atlases}\label{sec:atlas-construction}
\begin{wrapfigure}[12]{r}{0.5\textwidth}
\vspace{-\intextsep}
    \begin{minipage}{\linewidth}
    \begin{algorithm2e}[H]
        \footnotesize
        \setlength{\abovedisplayskip}{2.5pt}
        \setlength{\belowdisplayskip}{2.5pt}
        \setlength{\abovedisplayshortskip}{0pt}
        \setlength{\belowdisplayshortskip}{3pt}
        \caption{Atlas-DeePC}
        \label{alg:atlas-deepc}
        \LinesNotNumbered
        \DontPrintSemicolon
        \KwIn{%
            Offline data matrix $H$; feature index set $\calR$; 
            chart ranks $\{d_k\}_{k=1}^K$; 
            measured trajectory $w_{[t-L_\mathrm{sel},t-1]}$; 
            reference $y^{\mathrm{ref}} \in (\R^p)^{\Zp}$
        }
        1: \textbf{Construct the atlas:}
        Cluster $H(\calR,:)$ to obtain index sets $\{\calI_k\}_{k=1}^K$.
        Compute the SVD $H(:,\calI_k)=P_k\Sigma_kV_k^\top$
        and set $\Phi_k=P_k(:,1:d_k)$\;
        2: \textbf{Select a behavioral chart:}
        Observe the most recent trajectory
        $w_{[t-L_{\mathrm{sel}},t-1]}$
        and select $k_t$ using \eqref{eq:proj-res}\;

        3: \textbf{Apply the local policy:}
        \[
        u_t=\pi_{\calS_{k_t}}\!(
            w_{[t-\Tini,t-1]},
            y_{[t,t+\Tf-1]}^{\mathrm{ref}}
        )
        \]

        4: \textbf{Update:}
        Set $t\leftarrow t+1$\;

        5: \textbf{Repeat from line 2.}
    \end{algorithm2e}
    \end{minipage}
\end{wrapfigure}
This section develops Atlas-DeePC. Its modular structure consists of three steps:
(i) clustering offline trajectories, 
(ii) constructing the behavioral atlas by fitting a linear subspace to each cluster, and
(iii) online chart selection and predictive control. 
\paragraph{Behavioral atlas.}
Consider system \eqref{eq:ss} with restricted behavior $\calB_L$, where $L=\Tini + \Tf$. 
For each $j \in [1, N]$, let $w^{(j)}_{[0, N_j-1]}$ denote a trajectory collected from \eqref{eq:ss}. We refer to the $j^{\text{th}}$ data collection experiment as \textit{episode $j$}.
Let $T_j = N_j-L+1$. We construct the Hankel data matrix of episode $j$ as $H^{(j)} \coloneqq \calH_L\left(w^{(j)}_{[0, N_j-1]}\right) \in\R^{(m+p)L\times T_j}$. Here we use the permutation matrix $P$ to be consistent to DeePC \eqref{eq:deepc} row partition. 
Concatenating these gives the mosaic Hankel data matrix \citep{van2020willems}
\begin{equation}\label{eq:mosaic-hankel}
    H \coloneqq
    \left[
        H^{(1)}\ \cdots\ H^{(N)}
    \right]   =: 
    \left[
        h_1\ \cdots\ h_T
    \right]
    \in\R^{(m+p)L \times T}, \qquad 
    T\coloneqq\sum_{j=1}^N T_j.
\end{equation}
Let $\calR$ be an index set with $|\calR|=r$ and define $Z\coloneqq H(\calR,:) = [z_1\ \cdots\ z_T]\in\R^{r\times T}$. 
The vectors $z_i$ are called \emph{feature vectors} since they contain the relevant components of the input-output trajectories used for clustering (see Section~\ref{sec:reacher}).
In the next step, a data clustering method constructs a collection of $K>0$ nonempty sets of indices
\(
    \left\{
        \calI_1,\dots,\calI_K
    \right\}
\) based on the dataset $\{z_i\}_{i=1}^T$, 
where \(\calI_k\subseteq\{1,\dots,T\}\). 
For each $k=1,\dots,K$,
we define the local data matrix
\(H^{\calI_k} \coloneqq H(:,\calI_k)\).
The collection
\(
    \left\{
        H^{\calI_1},\dots,H^{\calI_K}
    \right\}
\)
therefore organizes the offline trajectories into local datasets.
We abuse the terminology of charts and atlases by analogy with differential geometry. 
A \emph{chart} provides Euclidean coordinates for a neighborhood on a manifold, and an \emph{atlas} is a collection of charts covering the manifold \citep{lee2012manifold}. 
We construct a low-dimensional linear subspace $\calS_k$ from $H^{\calI_k}$, which we call the $k^{\text{th}}$ \textit{behavioral chart}. For convenience, we also refer to $H^{\calI_k}$ as the $k^{\text{th}}$ behavioral chart when the distinction between the data matrix and its subspace representation is clear from context.
These charts form the \textit{behavioral atlas}
\(\left\{ \calS_1,\dots,\calS_K \right\}.
\)
For computation, we represent $\calS_k$ by an orthonormal basis matrix $\Phi_k$ such that $\im \Phi_k = \calS_k$.
Affine charts can also be constructed by fitting linear subspaces to centered local data and adding the corresponding offsets \citep{padoan2025affine}.
Atlas-DeePC is modular and accommodates various clustering methods, such as $K$-means \citep{mcqueen1967some} or $K$-subspaces \citep{bradley2000k}. 
We use $K$-means in the numerical examples in Section \ref{sec:reacher}, as it provides a simple implementation and yields strong performance. A systematic comparison of clustering methods is left for future work. 
Specifically, we partition the feature vectors $\{{z}_i\}_{i=1}^T \subset \R^r$ into $K$ clusters. The cluster \emph{label} for $z_i$ is denoted $\ell_i \in [1, K]$ and the \emph{centroid} of the $k^\text{th}$ cluster is denoted $c_k \in \R^r$. They are obtained by solving
\begin{equation}\label{eq:kmeans}
    \min_{\ell_1,\dots,\ell_T,\,c_1,\dots,c_K}
    \sum_{i=1}^T \| z_i-c_{\ell_i}\|_2^2.
\end{equation}
The resulting index sets $\calI_k\coloneqq\{i\in\{1,\dots,T\}:\ell_i=k\}$ form a disjoint partition of $\{1, \dots, T\}$. For computational efficiency, we approximately solve \eqref{eq:kmeans} using mini-batch $K$-means \citep{sculley2010web} with $K$-means++ initialization \citep{arthur2007k}. 
For each $H^{\calI_k}$, we compute its singular value decomposition
$H^{\calI_k}=P_k\Sigma_kV_k^\top$. 
Let $d_k \in [1, \rank H^{\calI_k}]$ be a prescribed rank, let
$\Phi_k\coloneqq P_k(:,1:d_k)$ and define the $k^\text{th}$ behavioral chart as $\calS_k = \im \Phi_k$. 
\paragraph{Online chart selection and predictive control.}\label{sec:online-control}
Let $L_{\mathrm{sel}} \leq L$ denote the number of recent samples used for chart selection.
At time step $t$, let $w_{[t-L_\mathrm{sel}, t-1]} \in \R^{(m+p)L_\mathrm{sel}}$ denote the most recent length-$L_\mathrm{sel}$ trajectory and define
\(\calS_k^{\mathrm{sel}} \coloneqq \mathrm{im}\left(\Phi_k(1:(m+p)L_{\mathrm{sel}},:)
    \right).\)
We select the chart with the smallest projection residual, i.e.,  
\begin{equation}\label{eq:proj-res}
k_t \in \underset{k\in\{1,\ldots,K\}}{\operatorname{argmin}}
\left\|
\left(I-\Pi_{\mathcal S_k^{\mathrm{sel}}}\right)w_{[t-L_{\mathrm{sel}},t-1]}
\right\|_2,    
\end{equation}
where $\displaystyle \Pi_{\calS_k^{\mathrm{sel}}}$ 
is the orthogonal projector onto $\calS_k^{\mathrm{sel}}$.  
If multiple charts attain the minimum residual, we select the chart with the smallest index. 
Given a reference output sequence
$y^{\mathrm{ref}}\in(\R^p)^{\Z_{\geq0}}$,
Atlas-DeePC (Algorithm \ref{alg:atlas-deepc}) yields the control policy
\begin{equation}\label{eq:atlas-policy}
u_t=\pi_{\calS_{k_t}}
\left(w_{[t-\Tini,t-1]},y^{\mathrm{ref}}_{[t,t+\Tf-1]}\right).
\end{equation}
\section{Simulation: Planar Robotic Manipulator}
\label{sec:reacher}
We evaluate Atlas-DeePC on a planar robotic arm tasked with moving its fingertip to a target position, with an additional F1TENTH \citep{o2020f1tenth} autonomous vehicle experiment reported in Appendix~\ref{sec:f1tenth}.
We first specify the atlas construction by $K$-means used in the experiments. 
All experiments were run on an Apple M1 Pro with an
8-core CPU and 16\,GB of unified memory.
\paragraph{Task.}
We use the Gymnasium \citep{towers2026gymnasium} \texttt{Reacher-v5} environment, a MuJoCo-based simulation  \citep{todorov2012mujoco}
of a two-link planar arm tasked with moving its fingertip to a target position, with sampling period $\Delta t = 0.02\,\text{s}$. 
Each simulation runs for $\Tsim=300$ steps and starts with an environment reset that randomizes the initial joint angles, joint velocities, and target position. 
The control input $u_t\in[-1,1]^2$ consists of torque commands applied to the two revolute joints. Let $q_{t}(1)$ and $q_{t}(2)$ denote the joint angles at time $t$, and let $p_t = (p_t^x, p_t^y) \in\R^2$ denote the position of the fingertip at time $t$. 
The default Gymnasium observation is defined as
\[
    y_t\coloneqq\bigl(\cos q_{t}(1),\cos q_{t}(2), \sin q_{t}(1),\sin q_{t}(2),\dot q_{t}(1),\dot q_{t}(2), p_t^x, p_t^y \bigr)\in\R^{8}.
\]
Given a target position $p^\mathrm{goal} \in \R^2$, the task is to drive the fingertip toward $p^\mathrm{goal}$ over
$\Tsim$ control steps.
Let $e_t\coloneqq\|p_t-p^{\mathrm{goal}}\|_2$ and
\(
\underset{t_1\le t\le t_2}{\operatorname{mean}} e_t
\coloneqq
\frac{1}{t_2-t_1+1}\sum_{t=t_1}^{t_2}e_t.
\)
We call a simulation successful if the mean fingertip error is below $10\,\mathrm{mm}$ over both the final 30 and final 100 time steps, and the error remains below 20 mm throughout the final 100 time steps, i.e., if
\[
\mean{T_{\mathrm{sim}}-\tau \le t\le T_{\mathrm{sim}}-1} e_t
<10\,\mathrm{mm} \quad \text{for } \tau \in\{30,100\}
\text{ and} \quad 
\max_{T_{\mathrm{sim}} - 100 \le t\le T_{\mathrm{sim}}-1} e_t
<20\,\mathrm{mm}.
\] 
The 100-step mean penalizes simulations in which the fingertip approaches the target only near the end of the simulation, while the 30-step mean checks accuracy near the end of the simulation. 

\begin{figure}[t]
    \centering
    \includegraphics[width=0.95\linewidth]{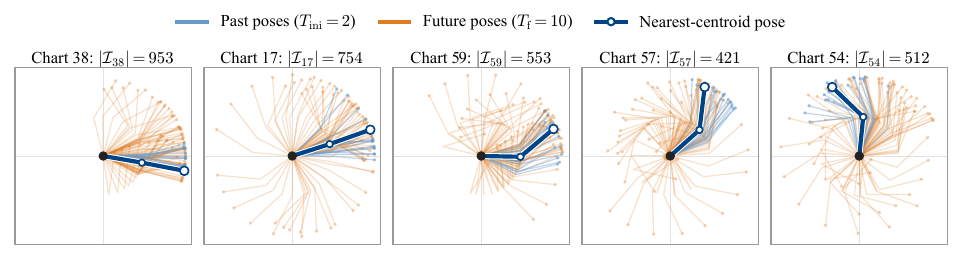}
    \caption{Sampled arm configurations from each of 5 behavioral charts. For visual clarity, 
    only a subset of the trajectories assigned to each chart is displayed. Each sampled trajectory contains $T_{\mathrm{ini}}=2$ past configurations (blue) and $T_f=10$ future configurations (orange).
    The bold blue arm shows the configuration corresponding to
    $h_{i_k^*}\bigl((m+p)(\Tini-1)+m+1:(m+p)\Tini\bigr)$,
where $i_k^* \in
\underset{i\in\calI_k}{\operatorname{argmin}}
\|z_i-c_k\|_2$,
i.e., the arm pose closest to the centroid in clustering coordinates.
}
\label{fig:physical_chart_clouds}
\end{figure}
\paragraph{Behavioral atlas.}
We collect trajectories from $N=200$ independent episodes, each of length $N_j = 200$. The inputs are sampled uniformly from $[-1,1]^2$. 
With $\Tini=2$, $\Tf=10$, $m=2$, and $p=8$, the resulting mosaic Hankel matrix \eqref{eq:mosaic-hankel} is $H\in\R^{120\times37800}$.
For atlas construction, we set the feature matrix $Z = H(\calR, :) \in \R^{4 \times T}$, where
\(
\calR \coloneqq \left\{
(m+p)(\Tini-1)+m+j
\;\middle|\;
j\in\{1,2,3,4\} \right\},
\)
i.e., the feature vectors consist of the first four components of the $\Tini$-th output of each trajectory $h_i$ in the Hankel matrix $H$. The feature vectors therefore specify a configuration of the arm.  
We cluster $Z$ using mini-batch $K$-means with $K=64$. 
For each cluster $k$, we set $d_k = d = 53$ and construct a behavioral chart $\calS_k = \im \Phi_k$ from $H^{\calI_k}$. 
The effects of chart rank $d$ and chart count $K$ are studied in Appendix \ref{app:rank-ablation}.
Figure~\ref{fig:physical_chart_clouds} visualizes the arm configurations associated with five charts.
\paragraph{Controller settings.}
The task specifies only the target position $p^\mathrm{goal} \in \R^2$, whereas the controller tracks an output reference in $\R^p$. We search for offline arm configurations with fingertip position closest to $p^\mathrm{goal}$, which is formulated as 
\[
(j^*, t^*) \in
\underset{
\substack{
j\in[1,N],\; t\in [0, N_j-1] \\}}
{\operatorname{argmin}}
\left\|
\left(y_t^{(j)}(7), y_t^{(j)}(8)\right)-
p^\mathrm{goal}
\right\|_2.
\]
We set $y^\mathrm{ref}_t = \bigl(
    y_{t^*}^{(j^*)}(1),\ldots,y_{t^*}^{(j^*)}(4),
    0,0,
    y_{t^*}^{(j^*)}(7),y_{t^*}^{(j^*)}(8)
    \bigr)$ for all $t \in [0, \Tsim-1]$. 
The cost weights in \eqref{eq:deepc} are
$Q=I_{\Tf}\otimes\diag(5\times10^3I_4,20I_2,10^5I_2)$
and $R=10^{-4}I_{m\Tf}$, with $\lambda_{\ini}=10^5$, where $\otimes$ denotes the Kronecker product. 
We impose the input constraint $\calU=[-1,1]^{m\Tf}$
and no explicit output constraints.
The quadratic programs are solved using Clarabel \citep{Clarabel_2024}.
We select the number of clusters $K$ and the chart rank $d$ by evaluating the control performance over 100 simulations and choosing the best-performing hyperparameters (see Appendix \ref{app:rank-ablation} for a more detailed analysis of the effects on the performance). The performance of Atlas-DeePC with the selected hyperparameters is then evaluated over 500 test simulations. 
\begin{figure}[t]
    \centering
    \includegraphics[width=0.9\linewidth]{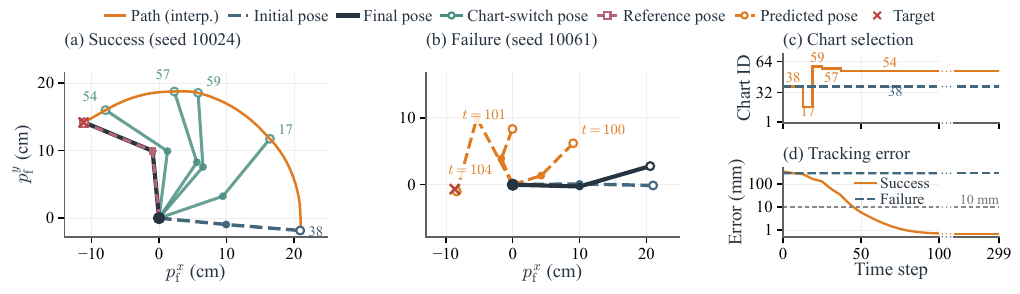}
    \caption{
    Atlas-DeePC simulations.
    (a) Fingertip path and arm poses at chart switches (green), labeled by the newly selected chart indices.
    (b) Failed simulation, with predicted poses (dashed orange line) from a single prediction horizon computed at $t=100$. Predicted joint angles and fingertip positions need not be kinematically consistent. The dashed poses show the predicted first-joint direction and fingertip position, without using the predicted second-joint angle.
    (c) Selected chart indices.
    (d) Fingertip errors $e_t$ on a logarithmic scale. 
    }
    \label{fig:reacher-baseline}
\end{figure}
\paragraph{Results.}
Figure \ref{fig:reacher-baseline} presents one successful and one failed simulation. 
In the successful case, Atlas-DeePC drives the fingertip to the target, with $\mean{T_{\mathrm{sim}}-30 \le t\le T_{\mathrm{sim}}-1} e_t < 1 \mathrm{mm}$ (Figure \ref{fig:reacher-baseline} (d)). 
In the failed case, the controller remains in one chart despite predicted poses approaching the target (Figure \ref{fig:reacher-baseline}(b)). This motivates a motion planning module that constructs intermediate references to guide the controller through behavioral charts toward the target.  
A statistical evaluation over 500 test simulations is presented in Section \ref{sec:reacher-continued}.

\section{Graph-Guided Reference Planning}\label{sec:reference-planning}
To guide the controller toward goals beyond the current operating region (i.e., towards trajectory references that are not in the selected behavioral charts), we select intermediate references using a weighted directed graph of observed chart transitions. 
\paragraph{Transition graph construction.}
We denote the set of successor pairs in the offline data by
\[
    \calT
    \coloneqq
    \left\{
        (i,i+1): i \in [1,T-1], 
        h_i \text{ and } h_{i+1}
        \text{ are columns of } H^{(j)}
        \text{ for some } j \in [1, N]
    \right\}.
\]
Using
the chart labels $\{\ell_i\}_{i=1}^T$ from \eqref{eq:kmeans}, define the transition count from chart $k$ to chart $k'$ as
\begin{equation}\label{eq:transition-count}
    n_{kk'}
    \coloneqq
    \left|\left\{
        (i,i+1)\in\calT
        \;\middle|\;
        \ell_i=k,\ \ell_{i+1}=k'
    \right\}\right|.
\end{equation}
We define the directed graph $\calG=(\calV,\calE)$ with node set
$\calV=\{1,\dots,K\}$ and edge set
$\calE\coloneqq\{(k,k'):k\neq k',\ n_{kk'}\geq n_{\min}\}$, where the user-defined threshold $n_{\min}\geq 1$ excludes rarely observed transitions. 
Each edge $(k, k') \in \calE$ is assigned the weight
$\omg_{kk'} \coloneqq
{d_\kappa(\calS_k,\calS_{k'})}/{\sqrt{n_{kk'}}}$,
where $d_{\kappa}$ is the chordal distance \citep{ye2016schubert, padoan2025distances} defined by 
$d_\kappa(\calS_k,\calS_{k'}) = \frac{1}{\sqrt{2}}\|\Phi_k \Phi_k^\top - \Phi_{k'}\Phi_{k'}^\top \|_F$, measuring the distance between subspaces $\calS_k$ and $\calS_{k'}$.
The edge weight is smaller for frequently observed transitions between similar behavioral subspaces.
\paragraph{Online reference planning.}
Let $k_t$ denote the chart index selected at time $t$.
Let $y_t^{\mathrm{term}} \in \R^p$ and $k_t^{\mathrm{term}}$ denote an application-dependent terminal reference and chart index.  
Let $\calP(k_t,k_t^{\mathrm{term}})$ denote the set of paths
from $k_t$ to $k_t^{\mathrm{term}}$ in $\calG$, where a \emph{path} $\rho$ is a finite ordered tuple of distinct nodes $\rho  = (v_0,\dots,v_M)$ where each $v_j \in \calV$ and $v_0=k_t$,
$v_M=k_t^{\mathrm{term}}$, and
$(v_{j-1},v_j)\in\calE$ for each $j\in[1,M]$. We denote the edge set of the path $\rho$ by $\calE(\rho) \coloneqq \{(v_{j-1},v_j):j\in[1,M]\}$.
Whenever $\calP(k_t,k_t^{\mathrm{term}})$ is nonempty, we compute a shortest path 
\begin{equation}\label{eq:shortest-path}
\rho_t^* \coloneqq 
    (k_t^0,k_t^1,\ldots,k_t^{M_t}) \in
    \underset{\rho \in \calP(k_t,k_t^{\mathrm{term}})}
    {\operatorname{argmin}}
    \sum_{(i,j)\in\calE(\rho)} \omg_{ij},
    \quad k_t^0 \coloneqq k_t, \quad 
    k_t^{M_t} \coloneqq k_t^\mathrm{term}
\end{equation} using Dijkstra's algorithm \citep{dijkstra1959note}.
The corresponding minimum path cost is denoted by 
\begin{equation}\label{eq:minimum-path-cost}
    d_\calG(k_t, k_t^\mathrm{term}) \coloneqq 
    \sum_{(i,j)\in\calE(\rho_t^*)}\omg_{ij},
\end{equation}
with cost zero when $k_t=k_t^{\mathrm{term}}$ and
$+\infty$ when no path exists. 
If $k_t\neq k_t^{\mathrm{term}}$, we select an intermediate output $y_t^{\mathrm{int}}\in\R^p$
from the next chart $H\left(:,\calI_{k_t^1}\right)$ in $\rho_t^*$ using
an application-dependent rule (see Section \ref{sec:reacher-continued}), and set
\(y^{\mathrm{ref}}_{[t, t+\Tf-1]}=
\mathbf{1}_{\Tf} \otimes y_t^{\mathrm{int}}\).
If $k_t=k_t^{\mathrm{term}}$, we set
\(y^{\mathrm{ref}}_{[t, t+\Tf-1]}=
\mathbf{1}_{\Tf} \otimes y_t^{\mathrm{term}}\).
As new measurements become available, we repeat the procedure using the newly selected chart $\calS_{k_{t+1}}$ from \eqref{eq:proj-res}, which need not coincide with the next chart $\calS_{k_t^1}$ on the shortest path. 
\section{Simulation: Reference Planning and Algorithm Comparison}
\label{sec:reacher-continued}
\paragraph{Reference selection for Reacher.}
We now describe how terminal and intermediate references are selected for the Reacher task introduced in Section \ref{sec:reacher}. We set $n_{\min}=3$ when constructing the chart-transition graph, so each edge corresponds to at least three observed transitions in the offline data.
Let $\bar{y}_i \coloneqq h_i\bigl((m+p)(\Tini-1)+m+1:(m+p)\Tini\bigr) \in \R^p$, which represents an offline arm pose. For each chart $\calS_k$, we compute and store 
\begin{equation}\label{eq:cluster-center-valid}
    i_k^*
    \in
    \underset{i\in\calI_k}{\operatorname{argmin}}
    \|z_i-c_k\|_2^2.    
\end{equation}
Thus, $\bar y_k^* \coloneqq \bar{y}_{i_k^*}$ is an observed arm configuration from cluster $k$ whose feature vector is closest to the centroid $c_k$. 
Since the centroid $c_k$ averages the feature vectors in cluster $k$ and need not correspond to an observed trajectory, we use $\bar{y}_{i_k^*}$ as the representative configuration of chart $\calS_k$ (Figure \ref{fig:physical_chart_clouds}).
Different arm poses can place the fingertip near the same target, but may differ in their reachability from the current chart. 
For a given target $p^\mathrm{goal}$, we choose ten offline outputs $\bar{y}_i$ whose fingertip positions $(\bar{y}_i(7), \bar{y}_i(8))$ are closest to the target $p^\mathrm{goal}$. 
Their indices form the candidate set $\calC$ defined by 
\[
    \calC
    \in
    \underset{
        \substack{\calJ\subseteq\{1,\ldots,T\}\\|\calJ|=10}
    }{\operatorname{argmin}}
    \sum_{i\in\calJ}
    \left\|
    \bigl(\bar y_i(7),\bar y_i(8)\bigr)
    -p^{\mathrm{goal}}
    \right\|_2.
\]
Each candidate indexed by $i \in \calC$ is associated with chart $\calS_{\ell_i}$, where $\ell_i$ is its cluster label from \eqref{eq:kmeans}. 
At each time $t$, we select the candidate whose chart has the minimum shortest-path cost from the current chart, i.e.,  
\[
    i_t^{\mathrm{term}}
    \in \underset{i\in\calC}{\operatorname{argmin}}\,
    d_\calG(k_t,\ell_i),
    \qquad
    k_t^{\mathrm{term}}\coloneqq\ell_{i_t^{\mathrm{term}}}.
\]
If multiple candidates attain the minimum path cost,
we select the one with the smallest fingertip distance to $p^{\mathrm{goal}}$. 
If these distances are also equal, we choose the smallest index.
The terminal reference is the selected output with its joint-velocity components set to zero, i.e., 
\[
    y_t^{\mathrm{term}}
    \coloneqq
    \bigl(
    \bar y_{i_t^{\mathrm{term}}}(1),\ldots,
    \bar y_{i_t^{\mathrm{term}}}(4),
    0,0,
    \bar y_{i_t^{\mathrm{term}}}(7),
    \bar y_{i_t^{\mathrm{term}}}(8)
    \bigr).
\] 
For the shortest path $(k_t^0,k_t^1,\ldots,k_t^{M_t})$ given by \eqref{eq:shortest-path}, we construct the intermediate output from the representative configuration of the next chart $\calS_{k_t^1}$
given by \eqref{eq:cluster-center-valid}:
\[
    y_t^{\mathrm{int}}
    \coloneqq
    \left(
        \bar y_{k_t^1}^*(1),\ldots,\bar y_{k_t^1}^*(4),
        0,0,
        \bar y_{k_t^1}^*(7),\bar y_{k_t^1}^*(8)
    \right).
\]
The reference supplied to Algorithm~\ref{alg:atlas-deepc}
is $y^{\mathrm{ref}}_{[t,t+\Tf-1]}
=\mathbf{1}_{\Tf}\otimes y_t^{\mathrm{term}}$
if $k_t=k_t^{\mathrm{term}}$ or
$\displaystyle \min_{0\leq\tau\leq t}e_\tau<70\,\mathrm{mm}$,
and $y^{\mathrm{ref}}_{[t,t+\Tf-1]}
=\mathbf{1}_{\Tf}\otimes y_t^{\mathrm{int}}$ otherwise.
We use a heuristic threshold of $70\,\mathrm{mm}$ to switch to direct terminal-reference tracking once the fingertip comes within this distance of the target, 
without claiming that this threshold is optimal.
Figure \ref{fig:recursive-term-ref} illustrates this procedure. At $t_a$, the next chart on the selected shortest path provides an intermediate reference (blue) toward the terminal pose (red). At $t_b$, the current chart has changed, leading to a different terminal candidate and an updated intermediate reference. At $t_c$, the current and terminal charts coincide, so the controller tracks the terminal reference directly.

\begin{table*}[t]
    \setlength{\tabw}{0.55\textwidth}
    \setlength{\colgap}{1.5em}
    \sbox{\reachertabbox}{%
        {%
        \footnotesize
        \setlength{\tabcolsep}{3pt}%
        \setlength{\fboxsep}{5pt}%
        \renewcommand{\arraystretch}{1.08}%
        \colorbox{tablebg!50}{%
        \begin{tabular}{@{}lrrrr@{}}
        \toprule
        Output dimension & \multicolumn{2}{c}{$p=8$} & \multicolumn{2}{c}{$p=10$} \\
        \cmidrule(lr){2-3}\cmidrule(lr){4-5}
        Reference Planning & No & Yes & No & Yes \\
        \midrule
        Chart count $K$ & 64 & 96 & 80 & 64 \\
        Selected rank & 53 & 57 & 64 & 62 \\
        Success rate & 90.0\% & 94.2\% & 83.8\% & \textbf{98.8\%} \\
        Median error (mm) & 1.81 & 1.86 & \textbf{1.51} & 1.92 \\
        90th percentile (mm) & 10.05 & 4.61 & 232.29 & \textbf{3.81} \\
        95th percentile (mm) & 218.27 & 12.12 & 296.36 & \textbf{4.63} \\
        99th percentile (mm) & 299.30 & 48.23 & 371.37 & \textbf{10.62} \\
        \midrule
        \multicolumn{5}{@{}l}{\textit{Median online computation time (ms/step)}} \\
        Chart selection & 0.16 & 0.21 & 0.22 & 0.19 \\
        Graph guidance & --- & 0.02 & --- & 0.02 \\
        QP solve & 19.81 & 19.65 & 22.77 & 22.41 \\
        Total & 20.01 & 19.92 & 23.03 & 22.65 \\
        \bottomrule
        \end{tabular}}}}%
    \centering
    \begin{minipage}[t]{\tabw}
        \vspace{0pt}
        \centering
        \usebox{\reachertabbox}
        \caption{Atlas-DeePC performance and computational cost on 500 test simulations comparing output dimensions and graph guidance. Error statistics summarize final-30 mean fingertip errors. 
        Chart selection measures the time spent computing projection residuals and selecting the behavioral chart.}
        \label{tab:reacher-heldout500}
    \end{minipage}\hspace{\colgap}%
    \begin{minipage}[t]{\dimexpr\textwidth-\tabw-\colgap\relax}
        \vspace{0pt}
        \centering
        \includegraphics
            [height=\dimexpr\ht\reachertabbox+\dp\reachertabbox\relax,
             width=\linewidth, keepaspectratio]
            {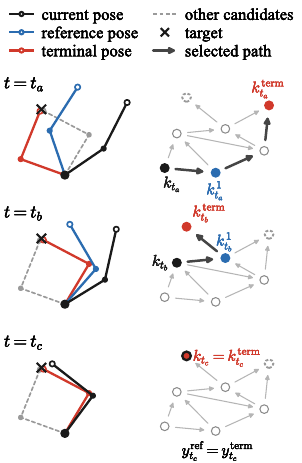}
        \captionof{figure}{Graph-guided reference selection. Intermediate (blue) and terminal (red) poses are updated as the current chart changes ($t_a \to t_b$).}
        \label{fig:recursive-term-ref}
    \end{minipage}
\end{table*}

\paragraph{Output augmentation.}
Alongside graph guidance, we study whether additional output coordinates improve tracking. 
Motivated by observable-based Koopman prediction and control
\citep{lian2021nonlinear, huang2023robust}, we augment the 8-dimensional output with the absolute orientation of the second link by
\(
    (
    y_t,\,
    \cos(q_t(1) + q_t(2)),\,
    \sin(q_t(1) + q_t(2))
    )\in\R^{10}.
\)
The augmentation is evaluated empirically.
Table \ref{tab:reacher-heldout500} compares Atlas-DeePC with and without graph guidance for output dimensions $p=8$ and $p=10$ over 500 test simulations, with chart counts and ranks selected separately based on development simulations reported in Appendix \ref{app:rank-ablation}. 
Here, ``error'' denotes
$\mean{\Tsim-30\le t\le\Tsim-1}e_t$.
For each fixed output dimension $p$, reference planning improves the success rate and percentile errors. 
Output augmentation alone does not improve performance, but its combination with reference planning achieves the best performance. 
These results suggest that adding output coordinates does not necessarily improve closed-loop performance. Appendix \ref{app:per-seed-results} isolates the effect of reference planning and provides additional individual simulation plots.
\begin{table*}[t]
    \noindent
    \begin{minipage}[t]{0.5\textwidth}
        \vspace{0pt}
        \raggedright
        \footnotesize
        \setlength{\tabcolsep}{4pt}
        \setlength{\fboxsep}{5pt}
        \colorbox{tablebg!50}{
        \begin{tabular}{@{}lcccc@{}}
            \toprule
            Method & \multicolumn{2}{c}{Atlas-DeePC}
                   & Select & Contextual \\
            \cmidrule(lr){2-3}
            Reference planning & No & Yes & No & No \\
            \midrule
            Coefficient dimension & 53 & 57 & 200 & 200 \\
            Success rate
                & 90.0\% & \textbf{94.2\%} & 83.0\% & 84.0\% \\
            Median error (mm)
            & \textbf{1.81} & 1.86 & 2.97 & 6.16 \\
            90th percentile (mm)
            & 10.05 & \textbf{4.61} & 29.95 & 11.66 \\
            95th percentile (mm)
            & 218.27 & \textbf{12.12} & 58.58 & 13.49 \\
            99th percentile (mm)
            & 299.30 & \textbf{48.23} & 83.48 & 91.92 \\
            \midrule
            \multicolumn{5}{@{}l}{%
                \textit{Total online computation (ms per step)}} \\
            Median
                & 20.01 & \textbf{19.92} & 1330.85 & 79.77 \\
            95th percentile
                & 27.36 & \textbf{23.15} & 3345.22 & 95.36 \\
            99th percentile
                & 53.85 & \textbf{43.17} & 4838.75 & 111.84 \\
            \bottomrule
        \end{tabular}
        }
    \end{minipage}\hfill
    \begin{minipage}[t]{0.42\textwidth}
        \vspace{0pt}
        \caption{Algorithm comparison on 500 test simulations
        (seeds 10001--10500). ``Select'' and ``Contextual'' denote select-DPC and contextual-DPC, respectively.
        Atlas-DeePC uses $p=8$, with $K=64$ and $d=53$ without
        guidance, and $K=96$ and $d=57$ with guidance.
        Coefficient dimension denotes the number of entries in the trajectory coefficient vector $g$ \eqref{eq:deepc}.
        Contextual-DPC uses the same fixed terminal
        references as unguided Atlas-DeePC.
        Error statistics summarize final-30 mean fingertip errors across simulations. 
        For Atlas-DeePC, the median online computation time is comparable to the $20\,\text{ms}$ sampling interval, although some steps exceed this interval. 
        }
        \label{tab:reacher-comparison}
    \end{minipage}
\end{table*}
\paragraph{Comparison with existing methods.}
We compare Atlas-DeePC with select-DPC \citep{naf2025choose} and contextual-DPC \citep{beerwerth2026less}. Both baselines select 200 trajectory columns from their respective offline datasets using relevance-based criteria and use these columns directly in the predictive control problem. Implementation details are provided in Appendix \ref{app:algorithm-comparison}. 
Our implementation of KDE-BDPC \citep{han2026big} did not achieve successful tracking on the Reacher benchmark and is therefore excluded from the main quantitative comparison. 
Appendix \ref{app:related-work} discusses the differences between its one-step control formulation and Atlas-DeePC. 
We conjecture that the one-step prediction of KDE-BDPC is not sufficient for long-term task planning required for the planar robotic task.  
Table \ref{tab:reacher-comparison} summarizes results over 500 test seeds. 
Similar to Table \ref{tab:reacher-heldout500}, ``error'' denotes $\mean{T_{\mathrm{sim}}-30 \le t\le T_{\mathrm{sim}}-1} e_t$.
In this benchmark, Atlas-DeePC achieves higher success rates and lower median fingertip errors than both Select- and Contextual-DPC, with substantially lower online computation times. 
Graph-guided reference planning further improves the success rate and reduces upper-percentile errors, while maintaining a similar computational cost. 
\section{Conclusion}
We presented Atlas-DeePC, a modular data-driven predictive control framework that constructs a behavioral atlas of low-dimensional charts from offline input-output trajectories. 
Online, the controller selects a chart from the behavioral atlas and applies the control policy \eqref{eq:atlas-policy}.  
An optional graph-guided reference planner uses observed chart transitions to generate intermediate references toward the target.
On the planar robotic manipulator benchmark, Atlas-DeePC achieved better performance and substantially lower online computation time than the compared methods. 
Future work will investigate principled choices of chart rank and chart count, approximation guarantees for nonlinear behavioral charts, and conditions for reliable closed-loop transitions between charts. We also plan to evaluate Atlas-DeePC on physical robotic platforms. 

\section*{Acknowledgments}
We used ChatGPT to assist with grammar checking, language editing, and code revision. 
The authors reviewed the resulting text and take responsibility for the final manuscript and implementation.
\bibliography{ref}
\newpage
\appendix
\section{Extended Related Work}
\label{app:related-work}
\paragraph{Behavioral systems theory and data-driven control}
Behavioral systems theory describes a dynamical system through its admissible trajectories, independently of a particular representation such as state space \citep{willems1986time, polderman1998introduction, willems2007behavioral}.
For controllable \gls{LTI} systems, the fundamental lemma \citep{willems2005note} gives conditions under which the set of all finite-length trajectories can be spanned by a Hankel matrix constructed from input-output data. 
This viewpoint provides a foundation for simulation and control directly from input-output data \citep{markovsky2008data, coulson2019data, markovsky2021behavioral}. Further developments address stabilization, optimality, robustness, and constraint satisfaction \citep{de2019formulas, berberich2021data, coulson2021distributionally, dorfler2022bridging, markovsky2022identifiability, sasfi2026soft}. 
Extensions to broader system classes include a fundamental lemma for linear parameter-varying systems \citep{verhoek2021fundamental}.
Complementary work on nonlinear and hybrid systems uses data partitioning and local affine-model fitting to identify piecewise-affine dynamics \citep{ferrari2003clustering, paoletti2007identification, bemporad2022piecewise}.
Nonlinear extensions also use basis functions and learned neural features to construct data-enabled predictors \citep{lazar2024basis, lazar2024neural}.
Further developments address computational efficiency through dimension reduction \citep{zhang2023dimension, shi2024efficient, shang2024convex, vahidi2025online, de2025scalable, alsalti_robust_2026}, online identification through subspace tracking \citep{sasfi2025great, jin2025online}, and online adaptation of data-based representations for nonlinear control \citep{teutsch2023online}.
\gls{DeePC} has been applied across a wide range of domains, such as power systems \citep{huang2021quadratic, huang2021decentralized},  and  robotics \citep{fawcett2022toward, kejun2024data, li2025hybrid, zeng2025reference, adang2026training, sharma2026rokon, wang2026velocity, ouyang2026direct, chen2026practical}.
\paragraph{Trajectory clustering for data-driven control}
Clustering provides a way to organize trajectories of a nonlinear system into local representations for prediction and control. 
Early work uses clustering to identify piecewise-affine models from input-output data \citep{ferrari2003clustering}. 
\citet{wu2014data} combine fuzzy operating-region memberships with subspace identification to construct predictive controllers. 
Gain scheduling constructs controllers for different operating conditions \citep{rugh2000research}. Following this approach,
gain-scheduled DeePC partitions the data into operating regions according to a prescribed scheduling variable \citep{guerrero2025gain, zieglmeier2025gain}. Choosing a suitable scheduling variable requires sufficient system knowledge to identify the dominant nonlinearity \citep[Assumption 2]{zieglmeier2025gain}.
In contrast, Atlas-DeePC clusters task-dependent features
of measured trajectories without requiring explicit
identification of the dominant nonlinearity.
More closely related big data-driven predictive control (BDPC) approaches cluster input-output trajectories into local sub-behaviors \citep{han2024big, han2025big, han2026big}. Their control formulations optimize a single future input-output pair conditioned on measured history \citep[Eq.~(11)]{han2024big}
\citep[Eq.~(45)]{han2025big}
\citep[Eq.~(34)]{han2026big}.
In contrast, Atlas-DeePC follows the standard multi-step receding-horizon formulation in predictive control. 
At each control step, it optimizes input-output trajectories over $\Tf$ future steps, applies the first optimal input, and repeats the optimization using updated measurements. 
Unlike one-step prediction, multi-step prediction allows the controller to evaluate control actions through their predicted effects over several future steps, reducing short-sightedness associated with optimizing only immediate tracking performance \citep{hatch2021value}.
\section{Additional Reacher Experiments}
\label{app:ablation}
\subsection{Chart counts and ranks}
\label{app:rank-ablation}
\begin{figure}[t]
    \centering
    \includegraphics[width=0.8\linewidth]{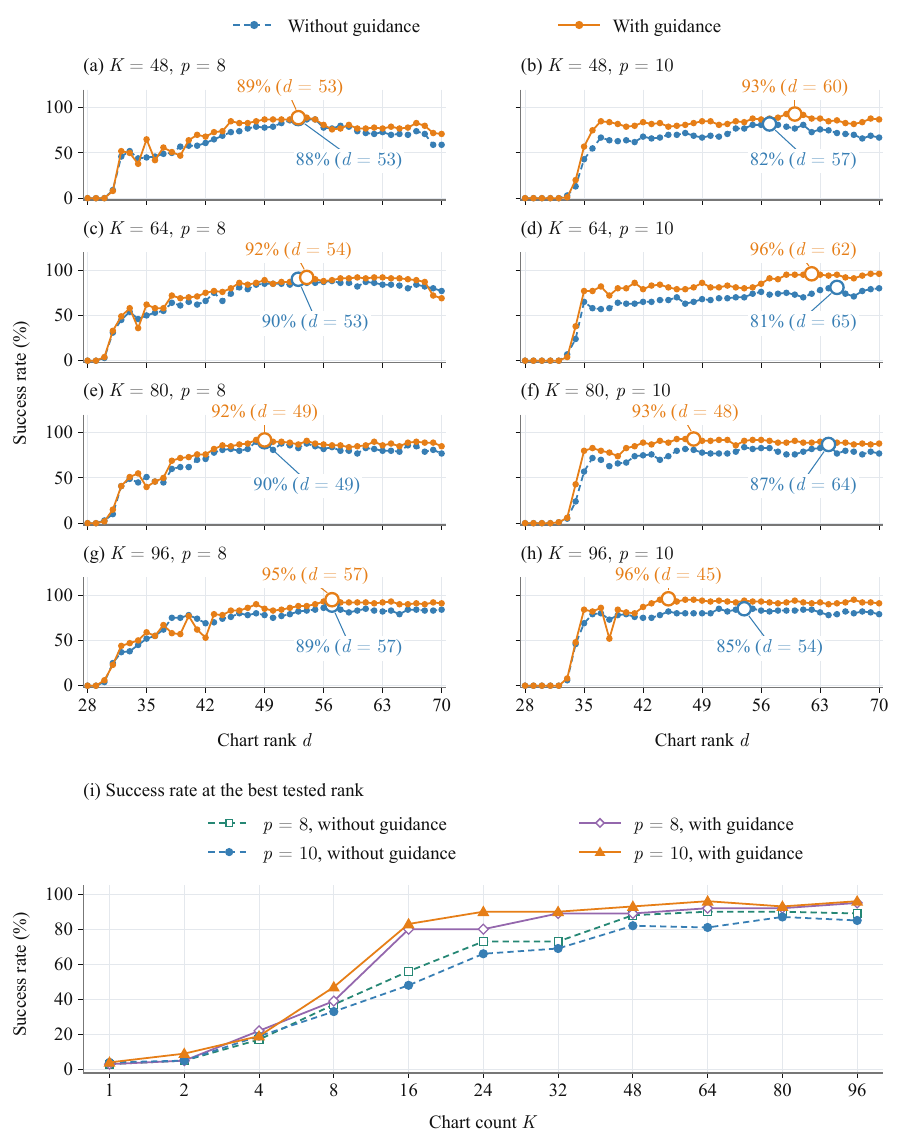}
    \caption{
    Success rate over 100 development seeds.
    (a) -- (h): Success rate versus chart rank $d$ for different chart counts $K$ and output dimensions $p$, with and without graph guidance. 
    Highlighted markers indicate the ranks selected by the highest success rate.
    (i): Success rate versus chart count $K$ for all four configurations, using the best rank at each $K$.
    }
    \label{fig:reacher-rank-ablation}
\end{figure}
The chart count $K$ determines how finely the offline trajectories are partitioned, while the chart rank determines the dimension of each local subspace. 
Figure \ref{fig:reacher-rank-ablation} reports 
success rates on 100 development simulations. 
Panels (a) -- (h) show success rates against chart rank for different chart counts and output dimensions, with and without graph guidance. 
Each panel shows success rates for chart ranks in $[28, 70]$. Lower ranks yield poor performance, whereas sufficiently large ranks produce a broad region of higher success rates. Increasing the rank further does not consistently improve performance.  
These development simulations inform the chart numbers and ranks used in Table \ref{tab:reacher-heldout500}. The unguided $p=8$ configuration uses $K=64$ and rank $d = 53$, as introduced in Section \ref{sec:reacher}. The other selected $(K, d)$ pairs are $(96, 57)$ for guided $p=8$, $(80, 64)$ for unguided $p=10$, and $(64, 62)$ for guided $p=10$. 
Table~\ref{tab:reacher-heldout500} evaluates these four
configurations on 500 test simulations, separately from the
development simulations used here for parameter selection.

\subsection{Effects of Reference Planning}
\label{app:per-seed-results}
\begin{figure}[t]
    \centering
    \includegraphics[width=\linewidth]{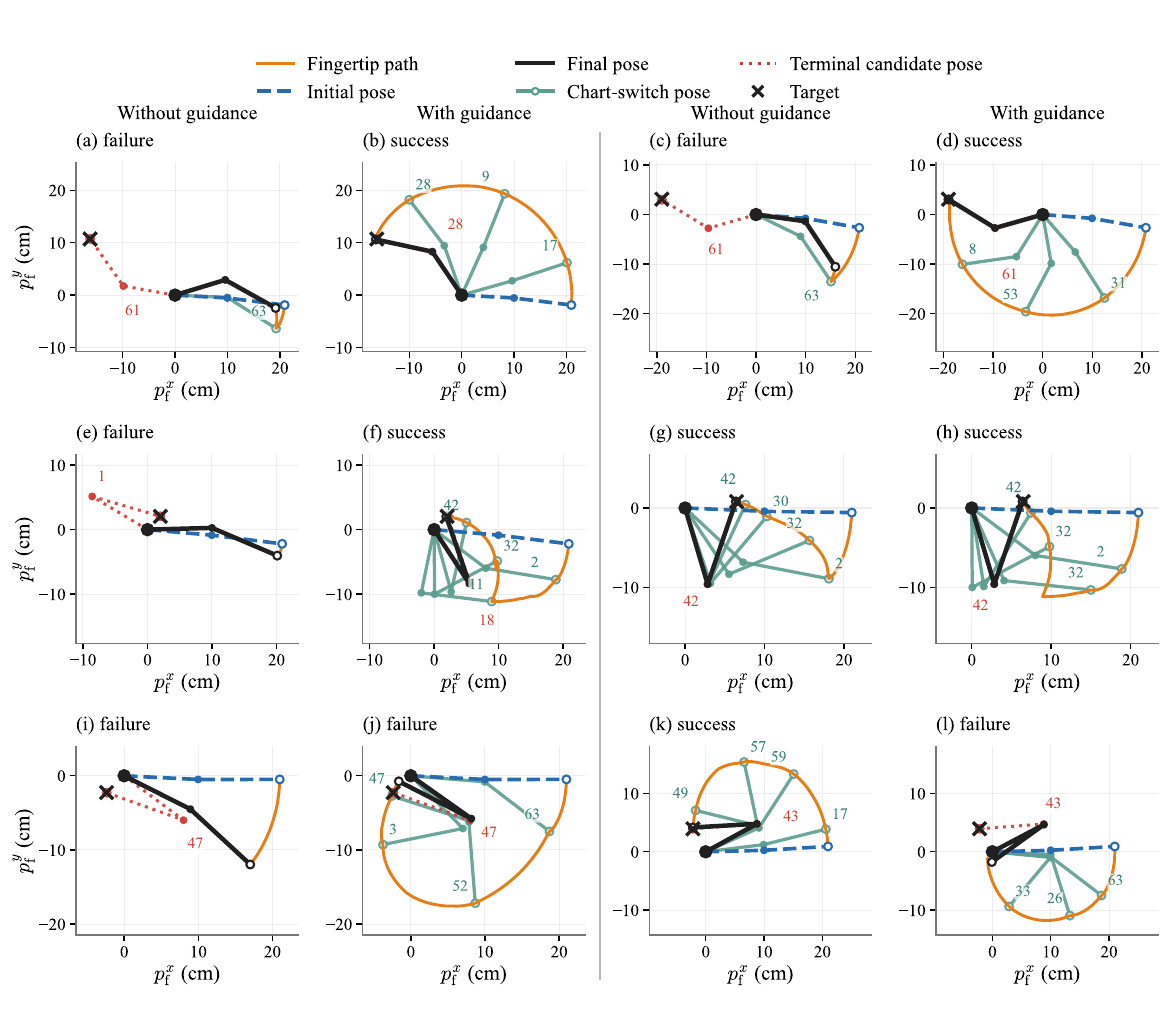}
    \caption{Effect of graph guidance on selected Reacher test simulations with $p=8, K=64, d=53$. 
    The first and third columns show simulations without guidance, whereas the second and fourth columns show simulations with guidance. 
    Numbers labeling the chart-switch and terminal reference poses denote their chart IDs. 
    These examples illustrate recovery, unchanged success or failure, and regression. }
    \label{fig:per-seed-plots}
\end{figure}

We isolate the effect of graph-guided reference planning by fixing $p=8, K=64, d=53$, and evaluating both configurations on 500 test simulations. 
Table \ref{tab:reacher-matched-guidance} shows that guidance increases the success rate from 90.0\% to 96.2\%, recovering 36 failed simulations while causing 5 previously successful simulations to fail. It also reduces upper-percentile errors. Figure \ref{fig:per-seed-plots} illustrates the effect of graph guidance on representative simulations. 
\begin{table*}[h]
    \centering
    \begingroup
    \footnotesize
    \setlength{\tabcolsep}{3pt}
    \begin{tabular*}{0.8\textwidth}{@{\extracolsep{\fill}}l*{6}{c}@{}}
        \toprule
        Graph guidance
        & \multicolumn{2}{c}{Success}
        & \multicolumn{4}{c}{Final-30 mean error (mm)} \\
        \cmidrule(lr){2-3}\cmidrule(l){4-7}
        &  & Rate & Median & $90^\text{th}$ percentile & 
        $95^\text{th}$ percentile & 
        $99^\text{th}$ percentile \\
        \midrule
        No & 450 & 90.0\% & \textbf{1.81} & 10.05 & 218.27 & 299.30 \\
        Yes & \textbf{481} & \textbf{96.2\%} & 2.04 & \textbf{4.49}
        & \textbf{7.40} & \textbf{38.00} \\
        \bottomrule
    \end{tabular*}
    \endgroup
    \caption{Effect of graph guidance on 500 test simulations (seeds 10001--10500), with $p=8$, $K=64$, $d=53$.
    Both configurations use the same atlas, controller parameters, initial states, and targets. Error statistics summarize final-30 mean fingertip errors across all simulations.}
    \label{tab:reacher-matched-guidance}
\end{table*}

\subsection{Implementation details of algorithm comparison}
\label{app:algorithm-comparison}
\subsubsection{Select-DPC and contextual-DPC}
Select-DPC \citep{naf2025choose} alternates between
data selection and predictive optimization at each
control step.  Our implementation uses the authors' publicly available code\footnote{\url{https://github.com/naefjo/choose-wisely-paper}}, adapted to our Reacher benchmark. 
In each iteration, it selects $n_s$
columns of the offline mosaic Hankel matrix $H$
based on the latest predicted input-output trajectory.
Denote the selected data matrix by $H_{\mathrm{sel}}$.
The optimization problem of Select-DPC is
\[
\begin{aligned}
    \min_{g,u^{\mathrm f},y^{\mathrm f},\sigma}\quad
    &(y^{\mathrm f}-y^{\mathrm{ref}})^\top Q
        (y^{\mathrm f}-y^{\mathrm{ref}})
     + (u^{\mathrm f})^\top R u^{\mathrm f}
     +\lambda_{\mathrm{ini}}\|\sigma\|_1
    +\lambda_1\|g\|_1
     +\lambda_\Pi\|(I-\Pi)g\|_2^2\\
    \text{s.t.}\quad
    &H_{\mathrm{sel}}g
    =
    \begin{bmatrix}
        u_{\mathrm{ini}}\\
        y_{\mathrm{ini}}+\sigma\\
        u^{\mathrm f}\\
        y^{\mathrm f}
    \end{bmatrix}, \quad 
    \mathbf 1^\top g=1,\quad
      u^{\mathrm f}\in\calU,\quad
      y^{\mathrm f}\in\calY.
\end{aligned}
\]
The projection regularizer \citep{markovsky2023data} is defined by
\[
    \Pi=M^\dagger M,
    \qquad
    M=
    \begin{bmatrix}
        H_{\mathrm{sel}}(\calJ,:)\\
        \mathbf 1^\top
    \end{bmatrix},
    \qquad
    \calJ=\{1,\ldots,(m+p)\Tini+m\Tf\}.
\]
Contextual-DPC \citep{beerwerth2026less} selects $n_s$
offline trajectories by measuring how closely their
past outputs match the recent measurements and their
future outputs match the reference.
The two distances are normalized and combined to rank
the candidate trajectories, see
\citep[Section~II-C and Algorithm~1]{beerwerth2026less}
for details.
Unlike Select-DPC, data selection and predictive
optimization are performed only once per control step.
The optimization problem of Contextual-DPC is
\[
\begin{aligned}
    \min_{g,u^{\mathrm f},y^{\mathrm f},\sigma}\quad
    &(y^{\mathrm f}-y^{\mathrm{ref}})^\top Q
        (y^{\mathrm f}-y^{\mathrm{ref}})
     + 
     (u^{\mathrm f})^\top R u^{\mathrm f}
     +
     (\Delta u^{\mathrm f})^\top R\,\Delta u^{\mathrm f}
     +\lambda_g\|g\|_2
     +\lambda_{\mathrm{ini}}\|\sigma\|_1\\
    \text{s.t.}\quad
    &H_{\mathrm{sel}}g=
    \begin{bmatrix}
        u_{\mathrm{ini}}\\
        y_{\mathrm{ini}}+\sigma\\
        u^{\mathrm f}\\
        y^{\mathrm f}
    \end{bmatrix},
    \qquad
    \mathbf 1^\top g=1,\quad
    u^{\mathrm f}\in\calU,\quad
    y^{\mathrm f}\in\calY,
\end{aligned}
\]
where
\(
    \Delta u^{\mathrm f}
    \coloneqq
    \bigl(
        u_t-u_{t-1},
        u_{t+1}-u_t,\ldots,
        u_{t+\Tf-1}-u_{t+\Tf-2}
    \bigr).
\)
\subsubsection{Implementation}
We evaluate both baselines on 500 test simulations. Implementation details are given below, 
and Table~\ref{tab:reacher-comparison} reports tracking performance and online computation times. 
\paragraph{Select-DPC.}
Since \citet{naf2025choose} include a Reacher benchmark, we use their publicly released implementation and dataset. 
We use the hyperparameters reported in
\citep[Appendix~D]{naf2025choose}:
\[
\begin{gathered}
Q=I_{\Tf}\otimes\operatorname{diag}(0,0,0,0,10,10,40000,40000),
\qquad R=10I_{m\Tf},\\
\lambda_1=10,\quad
\lambda_\Pi=10^4,\quad
\lambda_{\mathrm{ini}}=5\times10^6,\quad
\Tini=2,\quad \Tf=15,\quad
n_s=200,\quad n_{\max}=5,\\
\calU=[-0.5,0.5]^{m\Tf},
\qquad
\calY=
\left\{
    y^{\mathrm f}\in\R^{8\Tf}
    \;\middle|\;
    y^{\mathrm f}(8j+8)\leq0.1,\quad
    j=0,\ldots,\Tf-1
\right\}.
\end{gathered}
\]
Their dataset contains $200$ episodes of $200$
steps collected using i.i.d. inputs uniformly distributed
over $[-0.1,0.1]^2$. This dataset differs from ours, which uses exploration inputs in $[-1,1]^2$. 
We tested Select-DPC with our offline dataset, but did not achieve successful tracking. This is likely due to the hyperparameters in Select-DPC being data-dependent. 
We therefore used the authors' original dataset for comparison. 
Testing Atlas-DeePC on the original dataset of~\citet{naf2025choose} yielded similar results to those reported in Table~\ref{tab:reacher-comparison}. 
The trajectories are assembled into a mosaic-Hankel matrix
with $L=\Tini+\Tf=17$.
\paragraph{Contextual-DPC.} 
We adapt \citet{beerwerth2026less} to Reacher using our offline dataset.
Following their formulation, we stack the output and input as $(y_t, u_t) \in \R^{m+p}$ and apply the data normalization described in their Section II-D. 
We use \citep[Algorithm 1]{beerwerth2026less} for data selection. 
Following the output construction in
\citep[Section II-D, ``Incremental Input Representation'']{beerwerth2026less},
the tracking reference uses the same 
$y^{\mathrm{ref}}\in\R^8$ constructed in Section \ref{sec:reacher}.
The hyperparameters are tuned on 100 development simulations and set to  
\[
\begin{gathered}
Q = I_{\Tf}\otimes\operatorname{diag}
(0.01,0.01,0.01,0.01,10^{-4},10^{-4},100,100,0,0),
\quad R=0.1I_{m\Tf},\\
\lambda_g=20,\quad \lambda_{\mathrm{ini}}=100,\quad
\Tini=2,\quad \Tf=10,\quad n_s=200.
\end{gathered}
\]
For tracking the target, contextual-DPC uses the same fixed terminal reference as unguided Atlas-DeePC constructed in \mbox{Section~\ref{sec:reacher}}. \par

\section{F1TENTH Raceline Tracking}
\label{sec:f1tenth}  
This experiment evaluates Atlas-DeePC on the raceline-tracking task which is also considered by contextual-DPC \citep{beerwerth2026less}.
\begin{figure}[t]
    \centering
    \includegraphics[width=0.7\linewidth]
        {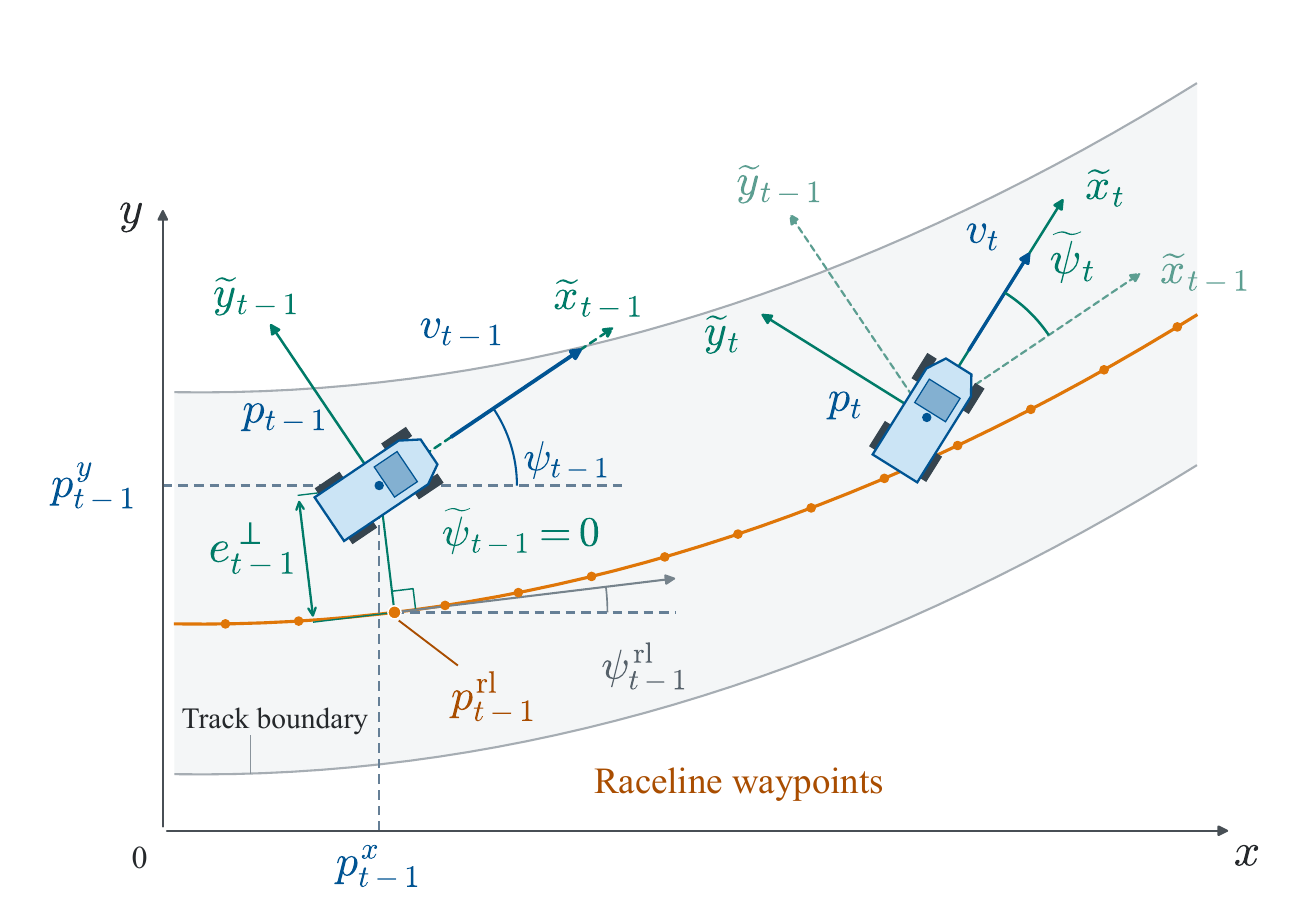}
    \caption{Vehicle motion and local coordinate frames.
    Each local frame is centered at the vehicle position
    and aligned with its heading. 
    At time $t$, the light-green
    dashed arrows show the directions of the previous frame,
    translated for comparison.
    Relative to that frame,
    $\widetilde{\psi}_{t-1}=0$ and
    $\widetilde{\psi}_t=\psi_t-\psi_{t-1}$.
    Orange dots denote raceline waypoints.
    The cross-track error $e_{t-1}^{\perp}$ is measured
    perpendicular to the raceline heading
    $\psi_{t-1}^{\mathrm{rl}}$ at the selected waypoint
    $p_{t-1}^{\mathrm{rl}}$.}
    \label{fig:f1tenth-tracking-schematic}
\end{figure}
\paragraph{Task.}
We consider raceline tracking in F1TENTH Gym \citep{o2020f1tenth} using nonlinear
single-track dynamics \citep{1024c5ad44034bbe9fae68e8f96ca260}
with sampling period $\Delta t=0.1\,\mathrm{s}$.
The input $u_t=(a_t,\delta_t)$ consists of acceleration $a_t$
($\mathrm{m/s}^2$) and steering angle $\delta_t$ (rad).
The output $y_t=(p_t^x,p_t^y,v_t,\psi_t)\in\R^4$ contains
the vehicle position $p_t=(p_t^x,p_t^y)$, speed, and heading,
as illustrated in Figure~\ref{fig:f1tenth-tracking-schematic}.
We express positions and headings in a local coordinate
frame. Let $\tau$ denote the time index of a sample and
$a$ the time index defining the frame, whose origin is
$p_a$ and whose positive $\widetilde x$-axis points along
$\psi_a$. The transformed position and heading are
\[
    \widetilde p_\tau =
    \begin{bmatrix}
        \cos\psi_a & \sin\psi_a\\
        -\sin\psi_a & \cos\psi_a
    \end{bmatrix}
    (p_\tau-p_a),
    \qquad
    \widetilde\psi_\tau=\psi_\tau-\psi_a.
\]
The components of $\widetilde p_\tau$ give the position
along the local axes, and $\widetilde\psi_\tau$ is the
heading relative to the positive $\widetilde x$-axis.
For each offline trajectory $h_i$, we transform all its
position and heading components using the position and
heading at the $\Tini$-th time step to define the local frame.
Online, we set $a=t-1$, using the latest measured position
$p_{t-1}$ and heading $\psi_{t-1}$ to define the frame
illustrated by the $\widetilde x_{t-1}$- and
$\widetilde y_{t-1}$-axes in
Figure~\ref{fig:f1tenth-tracking-schematic}.
The reference output sequence is constructed from the
raceline waypoints and prescribed speeds provided with
F1TENTH Gym. Both the recent measurements and the reference
output sequence are transformed using this same frame.
To measure tracking error, we select a closest raceline
waypoint, whose position $p_t^{\mathrm{rl}}$ and heading
$\psi_t^{\mathrm{rl}}$ are illustrated in
Figure~\ref{fig:f1tenth-tracking-schematic}.
The signed lateral error and its root-mean-square error
(RMSE) over a simulation of length $T_{\mathrm{sim}}$ are
\[
    e_t^\perp
    \coloneqq
    \begin{bmatrix}
        -\sin\psi_t^{\mathrm{rl}} &
        \cos\psi_t^{\mathrm{rl}}
    \end{bmatrix}
    (p_t-p_t^{\mathrm{rl}}),
    \qquad
    \mathrm{RMSE}
    \coloneqq
    \sqrt{\frac{1}{T_{\mathrm{sim}}}
    \sum_{t=0}^{T_{\mathrm{sim}}-1}(e_t^\perp)^2}.
\]
The vector $(-\sin\psi_t^{\mathrm{rl}},\cos\psi_t^{\mathrm{rl}})$
has unit length and points perpendicular to the raceline
heading. Thus, $e_t^\perp$ measures the vehicle's signed
lateral displacement from the selected waypoint, and
RMSE summarizes its magnitude over the simulation.  
A simulation lap is \emph{successful} if the vehicle completes one lap with \(|e_t^\perp| \le 1 \ \mathrm m\) for all $t$.
\paragraph{Behavioral atlas.}
\begin{figure}[t]
    \centering
    \includegraphics[width=0.9\linewidth]{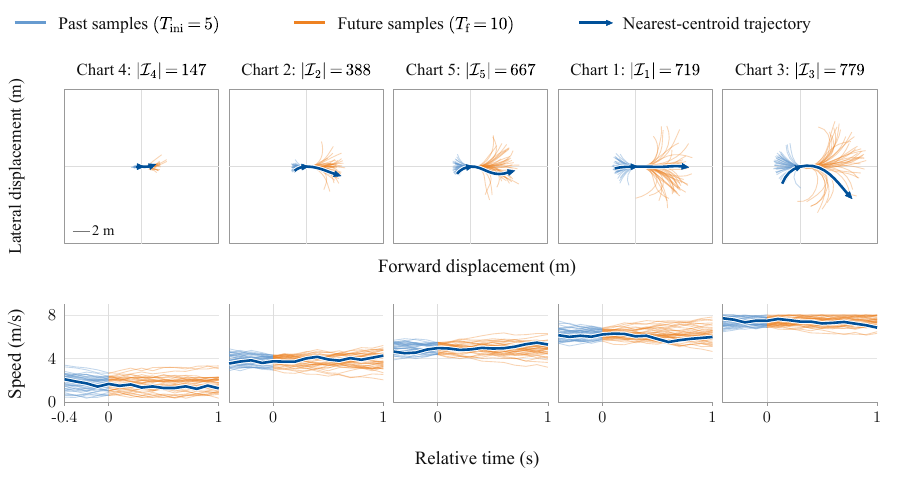}
    \caption{Behavioral charts of the noise-free F1TENTH data with $K=5$. Vehicle paths in local coordinates (top) and speed profiles (bottom). 
    Thin curves show $40$ sampled trajectories per chart, with initialization samples ($\Tini=5$) in blue and future samples ($\Tf=10$) in orange. 
    Bold dark-blue curves show the nearest-centroid trajectory $h_{i_k^*}$, 
    where $i_k^*\in\operatorname{argmin}_{i\in\calI_k}\|z_i-c_k\|_2$. Arrows indicate travel direction, and $|\calI_k|$ denotes the number of trajectories in chart $\calS_k$.}
    \label{fig:f1tenth_chart_clouds}
\end{figure}
We collect trajectories from $N=20$ independent episodes,
each of length $N_j=150$.
The initial speed is sampled uniformly from $[4,8]\,\mathrm{m/s}$.
Acceleration inputs are sampled independently and uniformly
from $[-3,3]\,\mathrm{m/s^2}$.
The steering inputs are generated by
\[
\delta_t=\min\{0.3,\max\{-0.3,\eta_t\}\},\text{ where }
    \eta_t=0.8\eta_{t-1}+0.6\xi_t
    \text{ and }
    \xi_t\overset{\mathrm{i.i.d.}}{\sim}\mathcal{N}(0,0.12^2).
\]
 Here, $\eta_{-1}=0$, and
$\xi_t\overset{\mathrm{i.i.d.}}{\sim}\mathcal{N}(0,0.12^2)$
means that $\xi_t$ are independent and identically
distributed Gaussian random variables with zero mean and standard deviation $0.12$.
At the sampling period of $0.1\,\text{s}$, independently sampled steering commands can cause abrupt steering changes and rapid changes in vehicle heading. We therefore generate temporally correlated steering inputs to obtain smoother steering variations. 
To assess robustness to sensor measurement noise, we add noise to the position, speed, and
heading measurements before coordinate alignment, yielding 
\[
    y_t=\bigl(
    p_t^x+\sigma s_1\epsilon_{t,1},\;
    p_t^y+\sigma s_2\epsilon_{t,2},\;
    v_t+\sigma s_3\epsilon_{t,3},\;
    \psi_t+\sigma s_4\epsilon_{t,4}
    \bigr),
    \qquad
    \epsilon_{t,j}\overset{\mathrm{i.i.d.}}{\sim}\mathcal{N}(0,1).
\]
Here, $\sigma \geq 0$ denote the noise level. 
The quantities $s_1, s_2, s_3, s_4$ are the standard deviations of $p_t^x, p_t^y, v_t, \psi_t$, respectively, computed across all time steps in the noise-free offline data before coordinate alignment.
We set $\Tini=5$, $\Tf=10$.
For atlas construction, the feature vector $z_i$ contains
only the speed at the $\Tini$-th sample of $h_i$,  i.e., 
\[
    \calR\coloneqq\{(m+p)(\Tini-1)+m+3\},
    \qquad Z=H(\calR,:)\in\R^{1\times T}.
\]
For each value of $\sigma$, we construct $Z$ from the noisy offline data and cluster it using mini-batch $K$-means,  
and construct a rank-$d$ behavioral chart
$\calS_k$ (Figure \ref{fig:f1tenth_chart_clouds}) from each local data matrix $H^{\calI_k}$ following the procedure in Section \ref{sec:atlas-construction}.
\paragraph{Controller.}
The controller tracks the prescribed raceline and speed
profile, expressed in the same local coordinate frame
as the measurements. The chart $\calS_{k_t}$ is selected by minimizing the projection residual \eqref{eq:proj-res} with $L_\mathrm{sel} = \Tini$.
The cost weights are
\[
    Q=I_{\Tf}\otimes\diag(100,100,1,1),
    \qquad
    R=I_{\Tf}\otimes\diag(0.1,10), \quad 
    \lam_\ini = 100.
\]
The constraints are
\(
    \calU=
    \bigl([-3,3]\times[-0.4189,0.4189]\bigr)^{\Tf},
\) and $\calY = \R^{4\Tf}$.
We use the control policy $\pi_{\calS_{k_t}}$ \eqref{eq:atlas-policy}, with the quadratic program solved using Clarabel.
\begin{figure*}[t]
    \centering
    \includegraphics[width=0.95\textwidth]{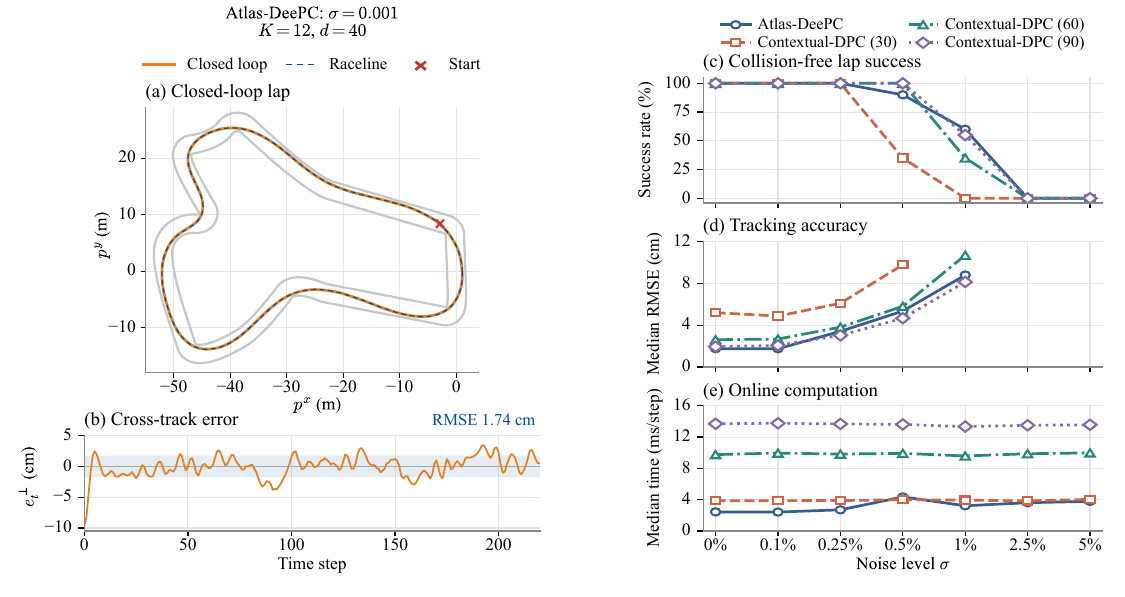}
    \caption{F1TENTH tracking and noise sensitivity.
    (a)--(b) Atlas-DeePC simulation at $\sigma=0.001$, with
    $K=12$, $d=40$, and seed 542: closed-loop path and signed cross-track error. 
    The shaded band indicates $\pm$ the simulation RMSE of $1.74\,\mathrm{cm}$.
    (c)--(e) Success rate, median cross-track RMSE over
    successful simulations, and median computation time over
    20 test simulations for each combination of noise level and controller configuration. 
    Missing RMSE points
    indicate no successful simulations. 
    Numbers in parentheses following Contextual-DPC indicate
    the number of data columns.}
    \label{fig:f1tenth-combined-results}
\end{figure*}
\paragraph{Results.}
We compare Atlas-DeePC with Contextual-DPC at noise levels
$\sigma\in\{0,\allowbreak 0.001,\allowbreak 0.0025,\allowbreak
0.005,\allowbreak 0.01,\allowbreak 0.025,\allowbreak 0.05\}$.
Contextual-DPC uses $30$, $60$, or $90$ selected trajectory
columns. At each noise level, both methods use the
same noisy offline dataset and are evaluated on 20 test
simulations. 
The Atlas-DeePC parameters $(K,d)$ are selected separately at each noise level using development simulations.  No reference planning is used.  
Figure~\ref{fig:f1tenth-combined-results} (a)--(b) shows an
Atlas-DeePC simulation at $\sigma=0.001$, with $K=12$, $d=40$. The vehicle completes a collision-free lap
with a cross-track RMSE of $1.74\,\mathrm{cm}$.
Panels~(c)--(e) summarize success rates, median cross-track
RMSE over successful simulations, and median online computation
times. Atlas-DeePC achieves tracking errors comparable to
Contextual-DPC using $90$ columns, while requiring
computation times generally lower than those of the $30$-column configuration.
\begin{table}[t]
    \centering
    \begingroup
    \footnotesize
    \setlength{\tabcolsep}{5pt}
    \begin{tabular}{@{}rrrccc@{}}
        \toprule
        $\sigma$ & $K$ & $d$
        & Maximum pairwise distance & Selected chart \\
        \midrule
        0      & 2  & 42 & $1.4142$             & 2 \\
        0.001  & 12 & 40 & $6.63\times10^{-15}$ & 1 \\
        0.0025 & 13 & 42 & $6.13\times10^{-15}$ & 1 \\
        0.005  & 14 & 60 & $7.16\times10^{-15}$ & 1 \\
        0.01   & 10 & 44 & $6.24\times10^{-15}$ & 1 \\
        0.025  & 3  & 68 & $4.73\times10^{-15}$ & 1 \\
        0.05   & 8  & 72 & $5.38\times10^{-15}$ & 1 \\
        \bottomrule
    \end{tabular}
    \endgroup
    \caption{F1TENTH chart-selection subspaces with
    $L_{\mathrm{sel}}=\Tini=5$ and $\Tf=10$.
    Maximum pairwise distance denotes
    $\max_{1\leq k<k'\leq K}
    d_\kappa(\calS_k^{\mathrm{sel}},\calS_{k'}^{\mathrm{sel}})$,
    computed using the chordal distance $d_\kappa$ defined in Section \ref{sec:reference-planning}. For each noise level, the column ``selected chart'' reports the chart index which remains constant over all time steps according to the chart selection rule~\eqref{eq:proj-res}. Chart 1 is selected at every time step in all 20 test simulations across all tested nonzero noise levels, whereas chart 2 is selected throughout all noise-free simulations.}
    \label{tab:f1tenth-selection-subspaces}
\end{table}
In these experiments, Atlas-DeePC selects the same chart throughout each simulation, i.e., chart 2 for $\sigma=0$ and chart 1 for all $\sigma > 0$ (Table \ref{tab:f1tenth-selection-subspaces}). 
For all tested nonzero noise levels, the selected subspaces $\calS_k^\mathrm{sel}$ are identical to numerical precision, therefore they yield the same projection residual in \eqref{eq:proj-res} for any trajectory $w_{[t-L_{\mathrm{sel}},t-1]}$, and the smallest-index
rule selects chart 1.  
This example demonstrates competitive tracking performance with low online computation, while showing that both methods remain sensitive to increasing noise.